\documentclass[12pt]{article}
\usepackage{graphics, latexsym, verbatim}
\newtheorem{thm}{\bf Theorem}[section]
\newtheorem{defi}[thm]{\bf Definition}
\newtheorem{prop}[thm]{\bf Proposition}
\newtheorem{lem}[thm]{\bf Lemma}
\newtheorem{rem}[thm]{\hspace{1em} \bf Remark}%[section]
\newtheorem{claim}[thm]{\hspace{3em} \bf Claim}%[section]
\newtheorem{cor}[thm]{\bf Corollary}%[section]
\newtheorem{exm}[thm]{\hspace{1em} \bf Example}

\newcommand{\bdry}{\partial}
\newcommand{\plane}{{\bf R}^2}
\newcommand{\rf}[1]{[{\bf #1}]}

\newcommand{\bpf}{{\bf Proof}\hspace{0.8em} }
\newcommand{\epf}{$\;\Box$}
\newcommand{\im}{{\rm Im}\,}
\renewcommand{\int}{{\rm Int}\,}
\newcommand{\inv}[1]{#1 ^{-1}}
\newcommand{\bR}{{\bf R}}
\newcommand{\bZ}{{\bf Z}}
\newcommand{\ser}[2]{#1_1,#1_2, \cdots , #1_#2}

\newcommand{\vph}{{\varphi}}
\newcommand{\ve}{{\varepsilon}}
\newcommand{\q}[2]{#1_1 ^2 + #1_2 ^2 + \cdots + #1_#2^2}
\newcommand{\brt}{{\bf R}^2}
\newcommand{\bdefi}{\begin{defi} }\newcommand{\edefi}{\end{defi}}
\newcommand{\blem}{\begin{lem} }\newcommand{\elem}{\end{lem}}
\newcommand{\bprop}{\begin{prop} }\newcommand{\eprop}{\end{prop}}
\newcommand{\brem}{\begin{rem} \rm  }\newcommand{\erem}{ \end{rem}}
\newcommand{\bclaim}{\begin{claim} }\newcommand{\eclaim}{\end{claim}}
\newcommand{\bexm}{\begin{exm} \rm  }\newcommand{\eexm}{\end{exm}}
\newcommand{\transv}{\,\top\hspace{-0.72em}{\raisebox{-0.1em}{$\cap$}}
\,\,}

\renewcommand{\marginpar}{\typeout}

\begin{document}

%\begin{comment}
\begin{flushleft}
{\large On the cusped fan in a planar portrait 
of a manifold
}

{\large Mahito Kobayashi}\\
Department of Information Science and Engineering\\
Akita University, 010-8502, Akita, Japan
\end{flushleft}
\vspace{2em}

{\small {\bf Abstract} 
A planar portrait of a manifold is the pair of the image and the
critical values of the manifold through a stable map into the plane. 
It can be considerd a geometric representation of the manifold drawn in
the plane. The cusped fan is its basic local configuration. 
In this article, we focus on the fibreing structure over the cusped fan,
and give its characterisation. 
As application, planar portraits of the real, complex, and quaternion
projective plane, regular toric surfaces, and some sphere bundles over
spheres etc. are constructed. 
Conversely, the source manifolds of certain planar portraits are detected. 
}
%\end{comment}

\section{Introduction}
\begin{figure}[h]\begin{center}
\scalebox{0.8}{\includegraphics{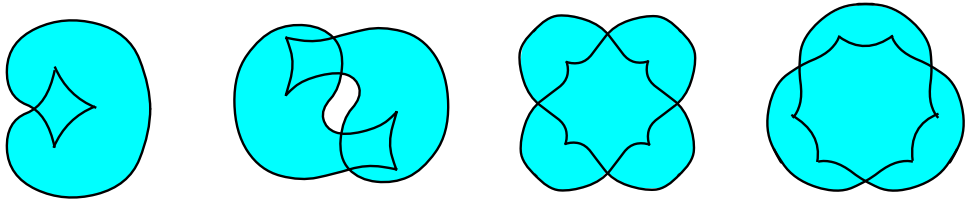}}
\caption{Planar portraits of certain closed manifolds}
\label{f:familiar}
\end{center}\end{figure}

How can one draw a picture of a closed manifold on the plane,
and what can one see from it ?
In this article, we make an attempt to draw and read pictures by using 
stable maps.
For a stable map $f:M \to \bR ^2$, $M$ a smooth closed manifold with dimension 
two or more, 
we call the pair 
$(f(M), f(S_f))$ up to diffeomorphism of $\bR^2$
a {\it planar portrait} of $M$ through $f$,
where $S_f$ denotes the set of singular points. 
Since the second factor $f(S_f)$ 
is a curve possibly with cuspidal points and normal 
crossings, the planar portraits have enough and not too much
complexity 
to be regarded a picture of $M$.
See Figure \ref{f:familiar} for samples of planar portraits.
We note that stable maps from closed $n$-manifolds ($n\geq 2$)
to 2-manifolds are known to be generic as smooth maps (Mather \rf{Ma}).

It is not straight to relate the planar portraits to manifolds.
In actual, 
infinitely many manifolds may admit a common planar portrait.
As a fundamental relation, 
a classical result by R. Thom \rf{Th} shows that the number of cusps in a
planar portrait is $\bmod \, 2$ congluent to the Euler characteristic 
$\chi (M)$ of $M$.
H. Levine \rf{L2} studied the indices of critical points 
of the composed Morse function $\gamma \circ f: M\to \bR^2 \to \bR$,
where $\gamma$ is a generic linear projection, viewed from the location
of $f(S_f)$.
R. Pignoni \rf{Pi} also obtained an Euler characteristic formula using 
$f(S_f)$ in case $M$ is a surface.
On the other hand, some examples 
show that finer topological properties are carried to the planar portraits
(see {\it e.g.}, 
applications to Corollary \ref{p:tuck}
in Section \ref{s:main}).
The purpose of this article is to give a step to read
views of maps and manifolds from the planar portraits,
and 
to offer samples of planar portraits as aids in 
seeking the suitable description for the relation
between the planar portraits and the manifolds.

\begin{figure}\begin{center}
\scalebox{0.75}{\includegraphics{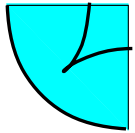}}
\caption{The cusped fan}
\label{f:cfan}
\end{center}\end{figure}

\begin{figure}\begin{center}
\scalebox{0.45}{\includegraphics{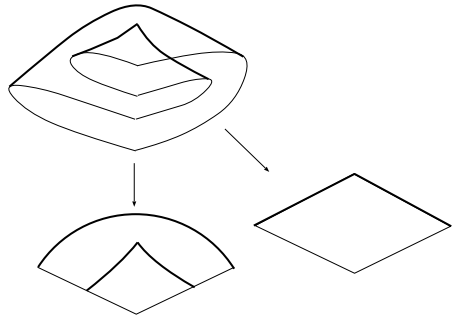}}
\caption{The cusped fan and the twice folding projection, in
	       dimension 2}
\label{f:2dimfold}
\end{center}\end{figure}

To read the planar portraits, one is required tasks of two kinds:
to know the fibreing structures over suitablly piecified planar portraits, 
and 
to know the glueings of manifold pieces over the subdivided 
planar portraits.
In this article, we deal with the first part and give a
characterisation of the fibreing structure over a basic piece 
named the {\it cusped fan} (See Figure \ref{f:cfan}, which is
strictly defined in Section \ref{s:main}).
In short, the fibreing structure is uniquely determined up to right-left equivalence
by the index of the cusp
and further it is a perturbation of a twice folding projection 
$(|z|^2,|w|^2)$ of $D^p\times D^q$, where $p$ and $q$ depend on the 
index of the cusp
(Theorem \ref{p:mainthm}. The twice folding projection 
is not stable in itself, as seen easily).
In the 2-dimensional case, what the result means
is easily understood: see Figure \ref{f:2dimfold}.

The result is applied effectively to both 
reading and drawing planar portraits. 
For the first, 
a tool for reading is given, and
admissible manifolds for 
certain planar portraits are detected
(Corollay \ref{p:tuck} and its application). 
For the second, lifting of Morse functions to stable maps
into $\brt$, and perturbations of orbit (quotient) maps of certain group actions
are considered
as tools for drawing.
In actual, planar portraits of projective planes $k P^2, k=\bR,
{\bf C}$ or ${\bf H}$, regular complex toric surfaces, and 
sphere bundles over spheres with cross-sections, are
constructed.
As a biproduct we can show that a certain planar portrait 
is realised by infinitely many maps, even if the source manifold
and the indices of singularities are fixed.
Some further applications 
are left to our proceeding papers  \rf{K2},\rf{K3}, {\it etc.},
for compactness.

The organization of this article is as follows.
In the next section, we introduce the cusped fan, states
the main theorem, and give its simple applications.
Through the proceeding four sections we prove the theorem, 
but technical calculations are left to the appendix.
In Section \ref{s:appli}, we give some applications of the theorem
to the construction of stable maps, and also give a biproduct
on the infiniteness of projections to a planar portrait.
One can skip to Section \ref{s:appli} after he read Section
\ref{s:main} to get the outline.
Henceforth,
manifolds and maps are assumed to be smooth
unless otherwisely mentioned.

\section{The cusped fan}\label{s:main}
A smooth map $f:M^n\to \bR^2 \, (n\geq 2)$ is called a {\it stable map}
if for any small perturbation $f'$ of $f$
in $C^\infty (M,\bR ^2 )$ equipped with Whitney  $C^\infty$ topology,
$f'=l\circ f \circ r$ holds 
for some self-diffeomorphisms $r$ and $l$ 
of $M$ and $\bR ^2$, respectively.
Two maps $f$ and $f'$ in this relation are called 
{\em right-left equivalent}.
For a stable map $f:M^n\to \bR^2 \, (n\geq 2)$, 
the set of singular points $S_f$ 
is made up of isolated
{\it cusps} and arcs of {\it folds} between the cusps.
It is known that 
in a small neighbourhood centered at a fold point
$P \in S_f$,  
the map $f$ is right-left equivalent to the 
normal form (F) bellow, 
where $a$ and $b$ are non-negative integers
(In case $a=0$ ({\it resp.}\, $b=0$), the terms $x$ and $|x|^2$
({\it resp.} $y$ and $|y|^2$) should be omitted).

\[(F):\quad
(u,x,y)\mapsto (u,|x|^2 - |y|^2),
u\in \bR, x\in \bR ^{a}, y\in \bR^{b}
\]

Similarly, a cusp point $P\in S_f$ has the normal form (C) bellow,
where $a$ and $b$ are non-negative integers
(In case $a=0$ ({\it resp.}$b=0$), the terms $y$ and $|y|^2$
({\it resp.} $z$ and $|z|^2$) should be omitted).

\[(C): \quad (u,x,y,z)\mapsto (u,x^3 -3ux+|y|^2 -|z|^2),
u,x\in \bR, y\in \bR ^a, z\in \bR ^b
\]

The {\em absolute index} of a cusp is the maximum of $\{a,b\}$
in $(C)$
(See \rf{L1} for the original, coordinate free definition).
As easily seen by $(F)$ and $(C)$,
the set of singular points $S_f$
is a $1$ dimensional proper submanifold of $M$, and 
the restriction of $f$ to $S_f$ except for the cusps
is an immersion.
Moreover the image of the immersion has 
at most normal crossings
and do not pass through the image of the cusps.
Refer to {\it e.g.}, \rf{GG} and \rf{L3} for basic notions on stable maps.

\begin{figure}[h]\begin{center}
\scalebox{0.75}{\includegraphics{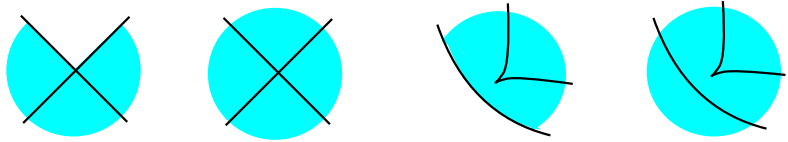}}
\caption{The duplicated local configurations in planar portraits}
\label{f:configs}
\end{center}\end{figure}

A {\it cusped fan} is 
the pair of a quadrant and a curve of two connected components
with one cuspidal singular point
pictured in Figure \ref{f:cfan}.
It can appear as a subpiece in 
a planar portrait $P=(f(M), f(S_f))$, where $f$ is a stable map,
as formed by a combination of a cusp with folds of definite types
$(u,|x|^2)$.
In constructing or in reading a planar portrait,
it is often convenient and useful to consider piecifications by
duplicated configurations 
(folds and folds, cusp and folds, cusp and cusp)
of critical values. 
These configurations are reduced to combinations of 
the four in Figure \ref{f:configs},
and the cusped fan is essentially the third case.
We note that fibreing structures for the first two
can be followed by using 
CPN (coordinatized product neighbourhood) in \rf{L3}
(where it is defined for dimension four).

\bdefi\label{d:scut} {\rm Let $X$ be a compact manifold with boundary, 
or possibly with corner.
A map $f: X \to \plane$ is 
a {\em fibrewise cut} of a stable map 
if there exist
an open manifold 
$\hat X$ 
which contains $X$ as a proper submanifold,
and a stable map 
$\hat f :\hat X \to \plane$
which have the three properties bellow:\\
1. The restriction $\hat f|X$ coincides with $f$,  \\
2. There exists a finite collection 
$\lambda_i, i=1,\cdots ,m$, 
of smooth plane arcs such that the intersections of $\lambda_i$ with other 
$\lambda_j$ and with 
$\hat f (S _{\hat f})$ 
are transverse, and\\
3. $X$ is obtained from $\hat X$ by cutting it along 
$\inv{\hat f}(\lambda _1)\cup \cdots \cup
\inv{\hat f}(\lambda _m)$.
}
\edefi

For a fibrewise cut of a stable map $f$,
we define the set of singularities
by $S_f = S_{\hat f}\cap X$, and 
call the pair $(f(X), f(S_f))$ the planar portrait of $X$, as before.

Our main result can be now stated as bellow:

\begin{thm}\label{p:mainthm}
Let $X$ be a compact connected $n$-dimensional manifold 
with corner, and 
$f : X \to \plane$ a fibrewise cut of a stable map
which has the planar portrait $(f(X), f(S_f))$ as
in Figure \ref{f:cfan}\,{\rm (}the cusped fan{\rm )}.
Then,
\begin{enumerate}
\item[{\rm 1}.] $X$ is diffeomorphic to $D^{a+1} \times D^{b+1}$ 
as a manifold with corner, where $a$ is the absolute index 
of the cusp and $b=n-2-a$.
\item[{\rm 2}.] $f$ is determined uniquely up to right-left equivalence,
by a given $a$.
\item[{\rm 3}.] We identify $X$ with $\{ |z|\leq 1, |w|\leq 1  \}$, 
where $z=(z_0, \cdots, z_a) \in \bR ^{a+1}$ and 
$w=(w_0, \cdots, w_b)\in \bR ^{b+1}$.
Then $f$ is right-left equivalent to the map 
\[(z,w) \mapsto (\,|z| ^2 +2\varepsilon\, k(z) w_0, |w| ^2 +2\varepsilon\, k(w) z_0\,)
\]
for any positive small constant $\varepsilon$, where $k(t)=1-|t|^2$.
\end{enumerate}
\end{thm}

\brem
Identify $X$ with $D^{a+1} \times D^{b+1}$ by $1$.
Then the boundaries 
$S^a \times D^{b+1}, D^{a+1} \times S^b$ 
and the corner 
$S^a \times S^b$ of $X$
are the inverse images of 
the two edges and the vertex of the quadrant $f(M)$, respectively,
by definition of a fibrewise cut, or by a direct calculation using 3.
\erem

Bellow is a simple application of the theorem. 
See section \ref{s:appli}, for more applications.

\begin{figure}[h]\begin{center}
\scalebox{0.9}{\includegraphics{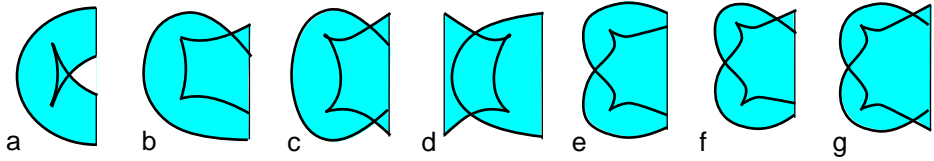}}
\caption{Bendings and tuckings in planar portraits}
\label{f:tucks}
\end{center}\end{figure}

\begin{cor}(Bending and Tucking lemma)\label{p:tuck}
Let $P=(h(X), h(S_h))$ be a planar portrait 
listed in Figure \ref{f:tucks}, where 
$h:X \to \bR^2$ 
is a fibrewise cut of a stable map. 
Then $X$ is diffeomorphic to 
the product $\Sigma^{n-1}\times D^1$ (in case of a,b,c,d),
$D^{n-1}\times S^1$ or  $D^{n-1}\tilde\times S^1$ (the non-trivial
 bundle) (in case of e),
or 
to $D^1\times S^1$ or $D^1\tilde\times S^1$ (in case of f,g),
where $\Sigma ^{n-1}$ is an $(n-1)$-homotopy sphere
in the 
gamma group $\Gamma _{n-1}$.
Conversely for each $P$,
every diffeomorphism types of $X$ above
can be realised.
\end{cor}

\begin{figure}[h]\begin{center}
\scalebox{0.75}{\includegraphics{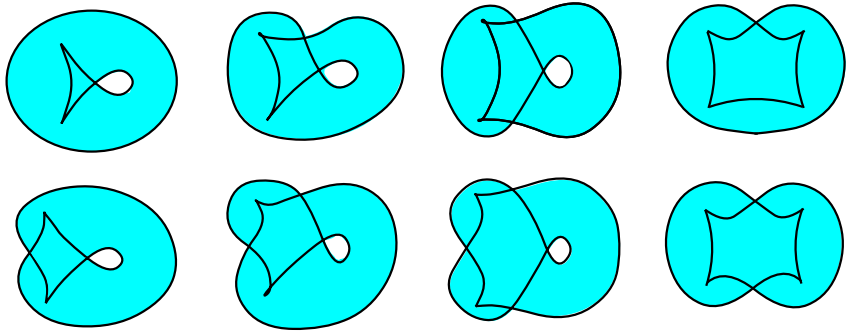}}
\caption{Planar portraits of $\Sigma ^{n-1}$ bundles over $S^1$}
\label{f:5ss1}
\end{center}\end{figure}

We can apply it to see, for example, that the source manifolds
of the planar portraits in  Figure \ref{f:5ss1} are all
diffeomorphic to the total space of a $\Sigma ^{n-1}$ bundle over $S^1$,
where $\Sigma^{n-1}$ is a homotopy sphere in $\Gamma _{n-1}$,
and conversely, each planar portrait is realised at least 
by $\Sigma ^{n-1} \times S^1$ 
and by the non-trivial bundle $S^{n-1}\tilde\times S^1$. 
The source manifold of the second planar portrait 
in Figure \ref{f:familiar}
is also diffeomorphic to a $\Sigma ^{n-1}$ bundle over $S^1$,
and the third one to either a torus or a Klein bottle.
Further the first one in  the figure 
is a planar portrait of $\bR P^2$ which is realised by 
combining the planar portrait e of a M\"{o}bius strip 
with a cusped fan.
In a similar way, various planar portraits can be constructed by 
combining 
the pieces listed in Figure \ref{f:tucks} and the cusped fans.

{\bf Proof of Corollary \ref{p:tuck}}
For the first three cases, one can move $h$ by homotopy 
to cancel the crossings of $h(S_h)$,
since these crossings are caused by folds of definite types,
and thus each fibre near the crossing point is disjoint.
Hence we can move $P$ to the union of two cusped fans
attached along an edge.
For the case d, we can move $P$ again to the same one as above,
by a half twist of the rectangular region enclosed by the cuspidal 
arcs and a side edge of $P$ by a homotopy of $h$.

Note that an arc of $S_h$ abuting the cusp consists of
folds of definite types. This implies that 
one edge of each cusped fan has inverse image diffeomorphic to
$D^{n-1}\times S^{0}$,
and hence $X$ is the union of two $D^{n-1}\times D^1$
attached along $S^{n-2}\times D^1$, by the theorem.
Namely,
$X$ is diffeomorphic to 
a $D^1$ bundle over $\Sigma$,
where $\Sigma$ is the union of two $D^{n-1}$ attached along boundary.
Its diffeomorphism class is either that of a trivial bundle
or the non-trivial one.
On the other hand, the restriction $h|\bdry X$ is 
a Morse function of $(n-1)$ dimensional manifold 
with four critical points of indices either $0$ or $n-1$,
which implies that 
$\bdry X$ is the disjoint union $\Sigma \coprod \Sigma$. 
Therefore $X$ is in actual the trivial bundle $\Sigma \times D^1$.

In the last three cases, we can move $P$ of f and g to e, 
similarly as before.
On the other hand, e is separated 
into the union of two degenerated cusped fans
along the cut line which passes through the crossing point.
Since the crossing is caused again by definite folds,
we can see that the inverse image of each piece above
is also $D^{n-1}\times D^1$.
Hence $X$ is the union of them
attached 
along $D^{n-1}\times S^0$,
or it is a total space of a $D^{n-1}$ bundle over $S^1$.
On the other hand, 
the normal form of the fold of definite types shows that 
the attaching diffeomorphism
in 
$\mbox{ Diff } (D^{n-1})$ preserves the norms
when $D^{n-1}$ is regarded a unit disc in $\bR ^{n-1}$,
and hence it can be linearlised, or moved to 
an $O(n-1)$ element by isotopy.
Hence $X$ has the required diffeomorphism types.
Note that $n=2$ in the case f and g, since both arcs abuting a cusp
are consisting of definite folds, and hence both boundary 
of the invese image of each cusped fan are diffeomorphic
to the disjoint union $D^{n-1}\coprod D^{n-1}$.
The construction part is easy (refer to Section \ref{s:appli}).
\epf

\section{A reduction to the smoothed case}\label{s:1ness}

\begin{figure}
\begin{center}
\scalebox{0.6}{\includegraphics{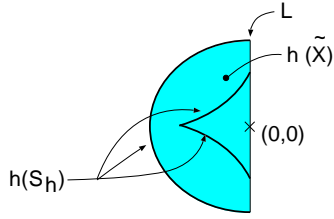}}
\end{center}

\caption{Planar portrait by the straightening $h$ of $f$}% in \ref{p:mainthm} }
\label{f:hdisc}
\end{figure}

We start to prove Theorem \ref{p:mainthm}.
To show assertion 1 and 2,
we will prove a corner removed version instead.
Suppose that $\tilde X$ is 
a compact connected $n$-manifold, $n\geq 2$ with boundary,
and $h:\tilde X \to \plane$ a fibrewise cut of a stable map 
with one cusp of absolute index $a$,
such that the planar portrait $(h( \tilde X ),h(S_h) )$ is as illustrated 
in Figure \ref{f:hdisc}{\rm :} $h( \tilde X )$ is the half-disc and
$h(S_h)$ is the bolded curves consisting of two components.
For simplicity, we assume that 
$h( \tilde X ) = \{ x^2+y^2 \leq 1, x \leq 0 \}$ and
set $L=\{0\}\times[-1,1]$.
Note that $\bdry \tilde X = \inv{h}(L)$, by the definition of a
fibrewise cut.

It is clear that
$\tilde X$ and $h$ can be obtained from 
$X$ and $f$ in Theorem \ref{p:mainthm}
by straightening of corner.
We can back to $f$ conversely, 
by restricting $h$ to the inverse image of a quadrant with vertex $(0,0)$.
The assertions are hence 
reduced to the proposition bellow.

\bprop \label{p:hdiscver}
Let $h: \tilde X \to \bR^2$ be a fibrewise cut of a stable map 
with its planar portrait as in Figure \ref{f:hdisc}. Then {\rm ;}
\begin{enumerate} 
\item[{\rm 1}.] $\tilde X$ is diffeomorphic to $D^n$, and 
$\bdry \tilde X$ is divided into two pieces 
$S^a \times D^{b+1}$ and $D^{a+1}\times S^b$ along 
the fibre $\inv{h}(0,0)$,
where $a$ is the absolute index of the cusp and $b=n-2-a$.
\item[{\rm 2}.] $h$ is determined uniquely up to right-left equivalence,
by $a$.
\end{enumerate}  
\eprop

In this section, we prove 1, and in the proceding section we do 2.
The first half of 1 is an easy application of \rf{L2}.
Hence we omit the proof.
To show the second half,
we regard $h|\bdry \tilde X$ 
as a Morse function onto $L$
by the lemma bellow.

\blem[Levine]\label{p:transMorse}
Let $f:M \to \plane$ be a stable map and
$\lambda$ an arc which is transverse to $f(S_f)$.
Then the restriction $f|\inv{f}(\lambda)$ can be considered
a Morse function onto $\lambda$, by composing an inclusion $\lambda 
\to \bR$.
\elem  

\bpf
This is just an
analogy of \rf{L3}, Proposition 2 in 1.3, in general dimensions. 
\epf

Note that $h|\bdry \tilde X: \bdry \tilde X \to L \subset \bR$ is 
a Morse function of the $(n-1)$-sphere with
four critical points.
By choosing the orientation of $L$ well, we may assume that 
their indices are
$0,a, a+1$ and $n-1$, lined in this order (\rf{L1}, Lemma(2) in 3.2),
and hence $\bdry \tilde X$ has a handlebody decomposition
$H_0 \cup H_a \cup H_{a+1}\cup H_{n-1}$.
We show here that on the middle regular level 
$\bdry (H_0 \cup H_a)$, the sphere 
$S'_a=\{0\}\times \bdry D^{n-1-a}$ 
in $H_a=D^a \times D^{n-1-a}$ (the 
attached sphere of the surgery of index $a$),
and the sphere 
$S_{a+1}=\bdry D^{a+1}\times\{0\}$
in $H_{a+1}=D^{a+1} \times D^{n-2-a}$
(the detached sphere 
of the surgery of index $a+1$)
meet transversely at one point.
Note that $H_0\cup H_a$ is a total space of an $S^a$ bundle 
over $D^{b+1}$ in general (recall that $n-1-a=b+1$).
But 
under the above transversality,
it is known to be a trivial bundle,
by standard argument in cancellation of handles
for $H_0\cup H_a \cup H_{a+1}$.
Hence $H_0\cup H_a$ is diffeomorphic to $S^a\times D^{b+1}$,
and $ H_{a+1}\cup H_{n-1}$ to
$D^{a+1}\times S^b$, by duality.

\begin{figure}
\begin{center}
\scalebox{0.8}{\includegraphics{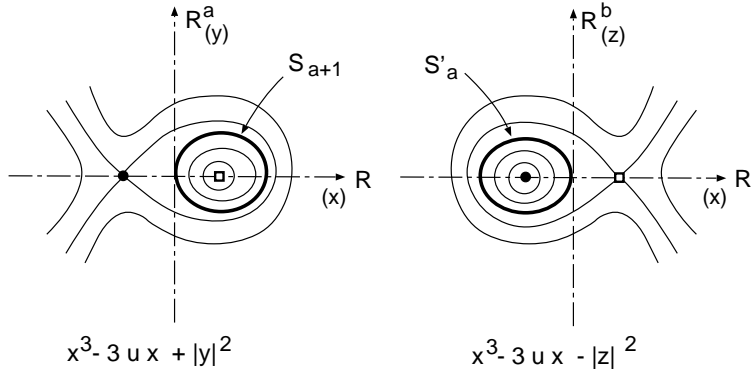}}
\end{center}
\caption{Attached and detached spheres in a neighbourhood of a cusp}
\label{f:cancelpair}
\end{figure}

To see the transversality of $S'_a$ and $S_{a+1}$, we may assume that $L$ 
is very close to the 
cusp image $h(P)$.
In actual, the Morse function $h|h^{-1}(L)$ 
changes only up to isotopy, when the line $L$
is moved close to $h(P)$ 
in keeping the transversality with $h(S_h)$.
Now by the normal form $(C)$, we can take a common neighbourhood
of the middle two critical points of $h|\bdry L$
and a system of coordinates, in which the Morse function is written by
$x^3-3ux+|y|^2-|z|^2$ for a positive small $u$.
It is easy to check that, by isotopy moves, 
the spheres $S'_a$ and $S_{a+1}$
can be placed in the hyper plane $|y|=0$ and $|z|=0$,
respectively, 
as shown in Figure \ref{f:cancelpair}.
In the figure, the dot and the box indicate 
the two critical points on the $x$ axis, 
one with index $a$
and the other with $a+1$, 
respectively.
Clearly these spheres meet transversely at the origin 
$(x,y,z)=(0,0,0)$.
\epf

\section{Uniqueness}
In this section, we prove $2$ of Proposition \ref{p:hdiscver}.
Denote by $P$ the unique cusp of $h$, and
take a small neighbourhood $U$ of $P$. 
We can take a system of coordinates
$(u,x,y,z)$ centered at $P$
such that on $U$,
$h$ takes the normal form $(C)$ in Section 
\ref{s:main}.

We modify CPN in \rf{L3} Proposition 1 in  2.2
by putting boundary conditions to it, as follows.

\begin{figure}
\begin{center}
\scalebox{0.65}{\includegraphics{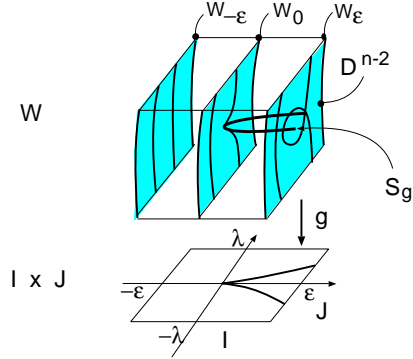}}
\end{center}
\caption{The cubic neighbourhood $W$}
\label{f:ccpn}
\end{figure}

\bprop\label{p:cubicnbhd}
Let $g(u,x,y,z), u,x \in \bR, y\in \bR ^a, z\in \bR^b$ 
be the normal form $(C)$ for the cusp.
For positive small numbers $\ve$ and $\lambda$,
there exist a product neighbourhood $I\times J$, 
$I=[-\ve ,\ve], J=[-\lambda ,\lambda]$
of the origin $0\in \bR^2$
and a
closed neighbourhood $W$ of the origin $0\in \bR^n$
diffeomorphic to 
$D^{n-2}\times I \times J$ 
such that the following properties hold 
{\rm (}refer to Figure \ref{f:ccpn}{\rm )}
{\rm :}
Let $W_u,  u\in I$ be the $u$-slice $W\cap (\{u\}\times \bR \times \bR ^a \times \bR ^b)$,
and  let
$g_u(x,y,z)$ be the function of $W_u$ given by $x^3 -3ux+|y|^2-|z|^2$. Then {\rm ;}

\begin{itemize}
\item[{\rm 1}.] $\{\ve \}\times J$ has two transverse intersections 
with the critical value set $g(S_g)$,
\item[{\rm 2}.] $g(W_u)= \{u\}\times J$,  
\item[{\rm 3}.] On the boundary of $W$ corresponding to 
$D^{n-2}\times I\times \{ \lambda \}$, each $g_u, u\in I$ is regular 
and takes constantly $\lambda$, 
\item[{\rm 4}.] On the boundary of $W$ corresponding to 
$D^{n-2}\times I\times \{ -\lambda \}$, each $g_u, u\in I$ is regular 
and takes constantly $-\lambda$, 
\item[{\rm 5}.] Each fibre $\inv{g_u}(t), u\in I , t \in J$ is transverse 
to the boundary of $W$ corresponding to 
$\bdry D^{n-2}\times I\times J$.
\end{itemize}
\eprop

We call $W$ a {\it cubic neighbourhood} of the cusp.
Its proof is similar to that in \rf{K1}, pp.347--348,
where the special case $a=b=1$ is treated,
and hence we omit it. \marginpar{bdry cond. 3,4 ?}
In the following, we regard $W$ a subneighbourhood of the cusp $P$ in 
the neighbourhood $U$.

\blem\label{p:nsuppt}
The connected component of $S_h$ which contains the cusp $P$ 
has a compact neighbourhood $N$ such as bellow {\rm :}
\begin{enumerate}
\item[{\rm 1}.] $N$ is diffeomorphic to the cubic neighbourhood
$W$
so that the boundary $N\cap\bdry\tilde  X$ is identified with the boundary $W_{\ve}$
of $W$, and
\item[{\rm 2}.] The restriction $h|N$ is right-left equivalent 
to the restriction $g|W$, where $g$ is the normal form of a cusp.
\end{enumerate}
\elem

\bpf
For a small neighbourhood of $P$ in $S_h$, 
the cubic neighbourhood $W$
itself deserves
as the neighbourhood $N$.
It is not difficult to enlarge such $N$ to the whole neighbourhood
of $S_h$, as seen bellow.
Take a linear projection 
$\gamma :\plane \to\bR$ 
so that  $\gamma \circ h$ is a Morse function 
with $\gamma\circ h (P)=0$ and with $\gamma (L)=1$
(refer to \rf{L2}).
Let $W_u$ be the $u$-slice of $W$ as in Proposition \ref{p:cubicnbhd}.
We may assume that $\gamma \circ h (W_u)=u$ for $u\in [-\ve,\ve ]$,
for a positive small number $\ve$.
Since 
the projection $h$ restricted to the right side of $\gamma ^{-1}(\ve)$
is simply  right-left equivalent to 
the product map $h|\inv{(\gamma\circ h)}(\ve)\times{\rm id}$,
where ${\rm id}$ is the identity map of $[\ve,1]$,
it is easy to 
enlarge $W$ 
to the required neighbourhood $N$ of the whole component of $S_h$.
\epf

Take the neighbourhood $N$ of $S_h$ as in Lemma \ref{p:nsuppt}, and 
we temporally replace $h$ 
inside $N$
with the projection 
corresponding to $D^{n-2}\times I \times J \to  I \times J$
of the cubic neighbourhood $W$.
Denote by $H:\tilde X \to \plane$ the resultant.
It is again a fibrewise cut of a stable map
but only with definite folds 
$(u,\ser{x}{{n-1}})\mapsto(u,\q{x}{{n-1}})$
as singular points.
It is clearly unique up to right-left equivalence,
since it can be shrinked preserving the right-left equivalence class, 
to a small neighbourhood of a definite fold point.
On the other hand,
$N$ has intersections with some regular fibres of $H$ %,
so that
it cuts $D^{n-2}$ from each fibre.
Hence any self-diffeomorphism of $\tilde X$
which preserves the fibres of $H$ is made into the identity map on a neighbourhood of $N$,
by moving it if necessary by a $H$-fibre preserving isotopy.

\begin{figure}
\begin{center}
\scalebox{0.7}{\includegraphics{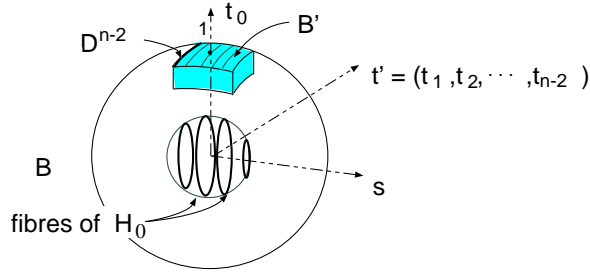}}
\end{center}
\caption{The neighbourhood $B'$ in the ball $B$}
\label{f:B}
\end{figure}

Denote by $N_e$ the closure of $\tilde X - N$.
The restriction $h|N_e$ is unique up to right-left equivalence.
In actual, 
$H$ is right-left equivalent to 
the map $H_0=(s^2 +|t|^2, s ):B \to \bR ^2$, 
where $B$ is the disc $s^2+|t|^2 \leq 1$ of $\bR ^n$,
and
at the same time, the restriction 
$H|N$ is identified with the restriction 
$H_0|B'$, where $B'$ is a small cubic collar neighbourhood
of the boundary point $s=0, t=(1,0,\cdots , 0)$
given by $(1-{\lambda}) ^2 \leq s^2+|t|^2 \leq 1, |s|\leq \ve $ and $|t'|\leq \mu  $
for positive small numbers $\ve , \lambda  $ and $\mu $,
where $t'=(t_1,t_2, \cdots , t_{n-2})$
(See Figure \ref{f:B}).
Hence $h|N_e$ is right-left eqivalent to the map
$H_0|\overline{B-B'}$.
Recall also that $h|N$ is uniquely determined up to right-left equivalence,
by a positive integer $a$ (2 of Lemma \ref{p:nsuppt}).
We now fix the integer $a$.

To show the uniqueness of $h$,
consider $h'$ which is another pasting of 
$h|N$ and $h|N_e$.
Then we replace $h'$ in $N$ as before and denote by $H'$ the resultant.
Note that 
$H'=\psi \circ H \circ \vph$ 
for self-diffeomorphisms $\psi$ of $\bR^2$
and $\vph$ of $\tilde X$.
But as mentioned before, $\vph$ is the identity map in a neighbouhrood of $N$.
Hence we have $h'=\psi \circ h \circ \vph$.
\epf

%%%%%%%%%%%%%%%
%%%%%%%%%%%%%%

\section{Jet transversalities}\label{s:transv}
Though this and the proceeding sections,
we prove $3$ of Theorem \ref{p:mainthm}.
Since the uniqueness property $2$ of the theorem is proved,
we are enough to check that the map presented in $3$ of the theorem %\ref{p:mainthm}
is a fibrewise cut of a stable map 
over a cusped fan with its absolute index of cusp $a$.
To save notation, we denote the polynomial described in $3$ of the theorem
also by $f$. 

Set $\Delta = \{ |z|\leq 1, |w|\leq 1  \}$. 
To show the stability of $f$,
we extend the domain slightly 
to an open neighbourhood $\hat\Delta$ 
of $\Delta$.
The stability of $f$ is equivalent
to the transversalities of jets 
$j^1 f \transv \Sigma _i$ 
and
$j^2 f \transv \Sigma _{ij}$ (see \rf{L1}),
besides with a certain genericity condition on 
the mutual position of the critical values.
In this section, we check the transversalities.
Only the outlines are stated here and the proofs
are given in Appendix.
Throughout the section, we set 
\begin{eqnarray*}
     p   &=&(1-2\ve w_0)z_0,  \\
     p_i &=&(1-2\ve w_0)z_i,  \, \mbox{ for } i=1,\cdots, a, \\
     q   &=& \ve k(z), \\
     r   &=& \ve k(w), \\
     s   &=&(1-2\ve z_0)w_0,  \mbox{ and } \\     
     s_j &=&(1-2\ve z_0)w_j,  \, \mbox{ for } j=1,\cdots, b.
\end{eqnarray*}

By straight calculation we can show:

\blem\label{p:singset}
Set $S=p\,s-q\,r$.
The set of singular points $S_f$ is defined by
$S=0, (\ser{z}{a})=0$, and $(\ser{w}{b})=0$.
\elem 
\marginpar{p:singset}

We can show further that (See Appendix for the proof):

\blem\label{p:Ssmooth}
$S_f$ is a regular curve in $\hat\Delta$
with two connected components.
\elem

Let $d^2 f : E \to L^* \otimes G$ be the second differential of $f$
defined in \rf{L1},
where $E = T\hat\Delta | S_1(f), \,G= {\rm cok}\, df$, and  $L =\ker df$.
The lemma bellow implies the transversality of 1-jet 
(3.1, \rf{L1}).
See Appendix for the proof.

\blem\label{p:1jet}
%$S_2 (f)= \emptyset$, 
At every point in $S_f$, ${\rm rk}\, df =1$ and $d\,{}^2 f$ is surjective.
\elem

Next we show 
the 2-jet transversality.
Denote by $S_{11}(f)$ the points of $S_f$
where the corank of $d\,{}^2f|L$ is $1$, by convention.
It is well known that a point $P\in S_{11}(f)$ is a cusp.

\blem\label{p:k=0}
Set $K=q\, \partial _{z_0}S-p \, \partial _{w_0}S$,
where $\partial _t S$ stands for $\displaystyle\frac{\partial S}{\partial t}$.
%for $t=z_0$ or $w_0$.
A singular point $P \in S_f$ is in $S_{11}(f)$ 
if and only if $K(P)=0$.
\elem

\blem\label{p:Ksingle}
$K$ restricted to $S_f$ vanishes only at a single point $P_D$ where $z_0=w_0$.
\elem

\blem\label{p:KStansvers}
At $P_D\in S_{11}(f)$, the locus $K=0$ is transverse to $S_f$.
\elem

One can see the proofs of these lemmaes in Appendix.
The 2-jet transversality follows the last one. % (refer to \rf{L1}).
For the absolute index of the cusp $P_D$, see Remark \ref{r:absidx} in Appendix.

\section{The planar portrait produced by $f$}\label{s:discrim}
The rest of proof is the check of 
the mutual position of the critical values
and the condition of a stable cut.
For a while until the global condition is checked, we continue to 
regard $f$ defined over 
an open neighbourhood $\hat\Delta$ of $\Delta$. 
Let $C_0$ and $C_1$ be the image of the component
of $S_f$ without, and with the cusp, respectively.

\blem\label{p:Csimple}
The two plane curves $C_0$ and $C_1$ are both simple and 
mutually disjoint.
\elem

\begin{figure}
\begin{center}
\scalebox{0.65}{\includegraphics{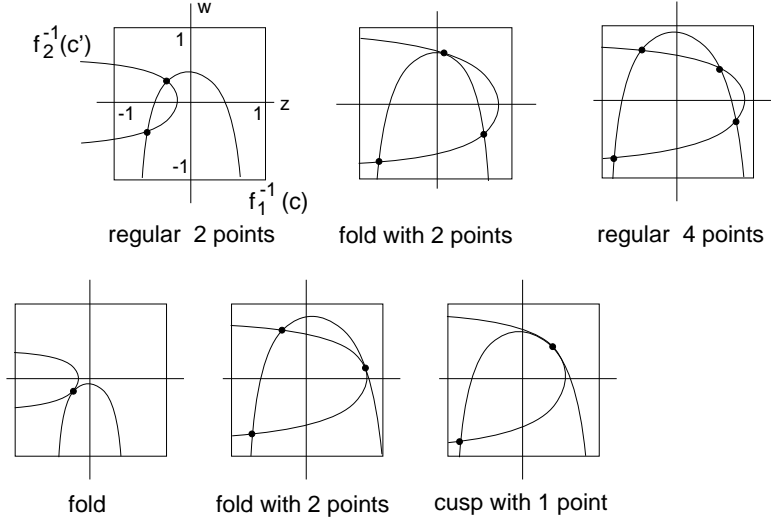}}
\end{center}
\caption{Fibres of $f$, in the 2-dimensional case}
\label{f:fibre}
\end{figure}

\bpf
By Lemma \ref{p:singset},
we can reduce the setup to the 
two variable case $z=z_0, w=w_0$.
Suppose that the curve $C_0 \cup C_1$ has a multiple point, 
say $f(Q)=f(Q')=(c,c')$ for some $Q,Q' \in S_f$.
Note that $c,c'<1$,
since $S_f \cap\bdry\Delta$ consists of the four points
$(\pm 1, 0)$ and $(0,\pm 1)$, 
which have mutually different images.
Write $f=(f_1,f_2)$.
Since 
$S_f$ is the locus of tangents of
the two level curves of $f_1$ and $f_2$,
the two curves $\inv{f_1}(c)$ and $\inv{f_2}(c')$  
have two or more tangent points.
But this is impossible
by the convexity of the curves (See Figure \ref{f:fibre}). 
\epf

By the lemma, we see that $f|S_f -\{\mbox{cusp}\}$ 
is an embedding and the cusp is not mapped onto
the image of the other singular points.
Hence the global condition of stability (see {\it eg.}, \rf{GG})
holds
and thus
$f:\hat\Delta \to \plane$ is stable.

We check 
that $f: \Delta \to \plane$ is a fibrewise cut of a stable
map.
Recall that $S_f$ is transverse to the boundary $\bdry\Delta$.
In actual,  $S_f$ lies in the $z_0 w_0$ plane (Lemma \ref{p:singset}),
and the normal direction $(\bdry_{z_0}S, \bdry_{w_0}S)$ to $S_f$ in the plane 
is not 
$(0,1)$ at $(0,\pm 1)$,
and is not $(1,0)$ at $(\pm 1,0)$,
as seen in the proof of Lemma \ref{p:Ssmooth} in Appendix.
Let $\lambda$ and $\lambda '$ be slightly extended arcs of
$\{1\}\times [-2\ve,1]$ and $[-2\ve,1]\times\{1\}$,
respectively.
The transversality of them to $f(S_f)$ follows 
from the above transversality and
the normal form $(F)$.
Then it is easy to see that
$f$ is a stable cut.

To see that $f$ produces a cusped fan, 
we are again enough to consider the two variable case (Lemma \ref{p:singset}). 
By the normal form $(C)$, 
the cusp image lies in the interior of $f(\Delta)$,
and 
hence 
$f(\Delta)$ is enclosed by $C_0, \lambda$ and $\lambda '$. 
We have two candidates for such planar portraits,
corresponding to the direction of the cusp: 
the cusp and nearby folds make an wedge, pointed to
$C_0$ or to the opposite.
Suppose that the second is the case.
Then we see a contradiction occurs by 
counting the number of inverse points $\inv{f}(Q)$ 
for a regular value $Q$:
it is two, if $Q$ is near $C_0$, 
inside the wedge, it should be greater by two, than that for
the outside by the normal form $(C)$, hence
it should be zero, outside the cuspidal edge.
It is a contradiction, 
since the corner point $(1,1)$ of the fan 
has four inverse points.

This ends the proof of Theorem \ref{p:mainthm}.

\section{Applications}\label{s:appli}
\begin{figure}\begin{center}
\scalebox{0.7}{\includegraphics{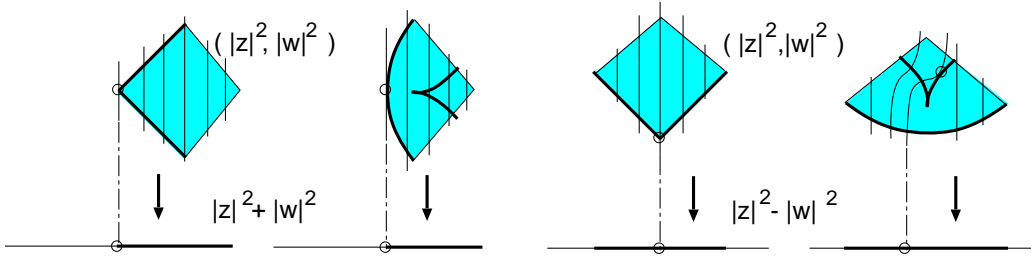}}
\caption{lifts of Morse critical points}
\label{f:lifts}
\end{center}\end{figure}

In this section, we give some applications of Theorem \ref{p:mainthm}.
The composition
$(z,w)\mapsto (|z|^2,|w|^2)\mapsto |z|^2\pm |w|^2$
shows that the twice folding projection is 
a lift into $\bR^2$, of a Morse type critical point of a function.
By applying the method in \rf{L2}, we see that
the stable cut over a cusped fan is still so, 
as illustrated in Figure \ref{f:lifts}:
The lines in the figure 
show the surjection
$\gamma :\bR^2\to \bR$ through which the Morse function is lifted. 
Let $f$ produce the cusped fan.
Then the tangent point of the line to $f(S_f)$ in the second figure
corresponds to the critical point of $\gamma \circ f$ of index $0$,
and that in the last figure corresponds
to the critical point of 
index $a+1$, where $a$ is the absolute
index of the cusp, by Theorem \ref{p:mainthm}.
By using this property, we can give 
lifts of some 
Morse functions.

\begin{figure}\begin{center}
\scalebox{0.7}{\includegraphics{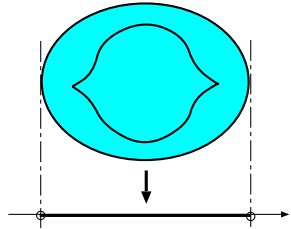}}
\caption{A lift of height function of the sphere}
\label{f:zwei}
\end{center}\end{figure}

\bexm\label{p:heightsphere}
Let $f$ be a stable cut over a cusped fan, of any absolute index of cusp.
We take two copies of $f$ and make a double to obtain a stable map 
$\tilde f: S^n \to \bR^2$, which produces the planar portrait as in Figure \ref{f:zwei}.
Clearly it is a lift of the height function of $S^n$
given by the projection of the unit sphere in $\bR^{n+1}$ to the first coordinate.
\epf
\eexm

\bexm\label{p:liftp2}
Let $g:kP^2 \to \bR, k=\bR, {\bf C},{\bf H}$ be a Morse function defined by
$g([x:y:z])=\displaystyle\frac{a|x|^2+b|y|^2+c|z|^2}{|x|^2+|y|^2+|z|^2}$,
where $0< a<b<c$.
It has three critical points at $[1:0:0], [0:1:0]$ and $[0:0:1]$,
of indices $0, l$ and $2l$ with critical values
$a,b$ and $c$, respectively,
where $l=1,2$ or $4$, corresponding to $k=\bR,{\bf C}$ or ${\bf H}$, respectively.

\begin{figure}\begin{center}
\scalebox{0.6}{\includegraphics{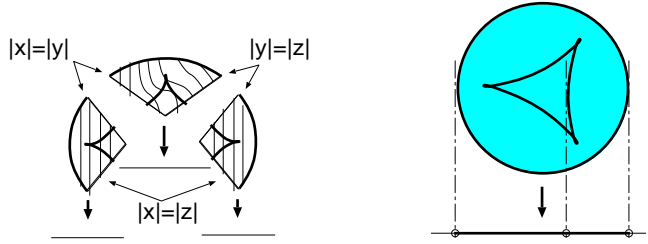}}
\caption{A planar portrait (right) of $kP^2, k=\bR,{\bf C}$ or ${\bf H}$}
\label{f:dritt}
\end{center}\end{figure}

Take a decomposition $kP^2 =H_0\cup H_1\cup H_2$,
where 
$H_0, H_1,  H_2$  
are neighbourhoods of 
$[1:0:0],[0:1:0]$ and $[0:0:1]$
defined by $\{
|y|,|z| \leq |x|\},
\{|z|,|x| \leq |y|
\}$
and $\{|x|,|y| \leq |z|\}$, respectively.
Identify each of them with
$\{ |s_i|, |t_i|\leq 1\}$ and hence to $D^l\times D^l$,
where 
$s_0=\displaystyle\frac{y}{x}, t_0=\frac{z}{x}$,
$s_1=\displaystyle\frac{z}{y}, t_1=\frac{x}{y}$,
and
$s_2=\displaystyle\frac{x}{z}, t_2=\frac{y}{z}$.
On these pieces, $g$ is right equivalent to $s_0 ^2 + t_0 ^2, s_1 ^2 - t_1 ^2$,
and $-s_2 ^2 -t_2 ^2$, respectively.
Take three copies of $f:D^l\times D^l \to \bR^2$ in 3 of 
Theorem \ref{p:mainthm}, and place the cusped fans on the plane 
as in Figure \ref{f:dritt}, left,
in refering to Figure \ref{f:lifts}.
It is easily seen from normal form 3 in Theorem \ref{p:mainthm}
that these map pieces can be well pasted to 
form a map $\tilde f: kP^2 \to \bR^2$ which produce the planar portrait
as in Figure \ref{f:dritt}, right, by moving the images up to isotopies, if necessary.
It is clearly a lift of $g$, as shown in the figure. \epf
\eexm

The stable cut over the cusped fan 
can be also regarded a perturbation of the quotient map 
of the linear product actions of $O(p)\oplus O(q)$  
to $D^p \times D^q$.
By using this, we can obtain  
a perturbation of the moment map of a regular toric complex surface
and the planar portrait produced by them.

\bexm\label{p:toric}
For $\gamma \in \pi _1 SO(2)=\bZ$,
denote by 
$\tilde \gamma$ the diffeomorphism of
$D^2 \times \bdry D^2 \to \bdry D^2 \times D^2$
defined by
$\tilde \gamma(z,w)=(\bar w, w ^{\gamma} z)$,
where
$D^2$ is regarded as a unit disc in the $\bf C$ plane,
and 
$\bar w$ denotes the complex conjugate $(w_0,-w_1)$
of $w=(w_0,w_1)$.
It is well-known that a regular toric surface $M$ is obtained by 
putting $k$-copies of $D^2\times D^2$ together,
where $k$ is the Euler characteristic $\chi (M)$ of $M$,
by using the diffeomorphisms $\tilde \gamma _i (z,w)=
(\bar w, w ^{\gamma _i} z), i=1,\cdots , k$
for some $\gamma_1, \cdots , \gamma _k \in \bZ$. 

\begin{figure}\begin{center}
\scalebox{0.6}{\includegraphics{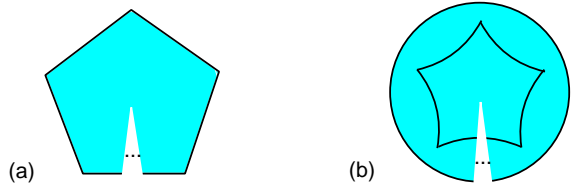}}
\caption{(a) image of the moment map, and (b) the planar portrait
of a regular toric surface}
\label{f:asa}
\end{center}\end{figure}

Take $k$ copies of $f$ 
given in $3$ of Theorem \ref{p:mainthm}.
By composing certain isotopy of $\bR ^2$ to $f$,
we may assume that $\im f$ is a sector with
vertex of angle $2\pi /k$.
Then as before in Example \ref{p:liftp2},
one can obtain a stable map $\tilde f:M \to \bR^2$
with $k$ cusps such that the planar portrait is obtained by pasting
$k$ cusped fans to form a disc
(Figure \ref{f:asa}, (b) ).

We note that $\tilde f$ is a perturbation of a moment map:
Let $q: M\to P_k\subset \bR^2$ be the map obtained by 
taking the unfolding parameter $\ve$ in $f$ to be zero
in the above construction,
where $P_k$ is a polygon with $k$ edges (Figure \ref{f:asa},(a)).
It is clearly a moment map, or the quotient map
of the toric action. 
\epf
\eexm

\brem
For the projective plane $kP^2$, the maps in Example
\ref{p:liftp2} and in Example \ref{p:toric} are the same.
Note that the Morse function $g$ in Example \ref{p:liftp2}
can be factorized as
\[
 [x:y:z]\mapsto \frac{1}{|x|^2+|y|^2+|z|^2}(a|x|^2, b|y|^2, c|z|^2)
 \mapsto  g([x:y:z]),
\]
where 
the first factor is a moment map of  $kP^2$ 
with the image a triangle in $\bR ^3$, which we regard a map into 
$\bR^2$.
The map in Example \ref{p:toric} is the perturbation of it.
\erem

Bellow is a biproduct of this construction.

\begin{figure}\begin{center}
\scalebox{0.6}{\includegraphics{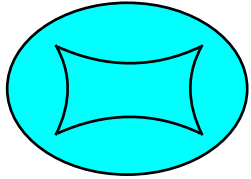}}
\caption{A common planar portrait of 
sphere bundles over spheres with cross sections}
%Hirzebruch surfaces}
\label{f:quatro}
\end{center}\end{figure}

\begin{cor}\label{p:infinite}
Let $M$ be either $S^2\times S^2$ or $S^2\tilde\times S^2$,
where the latter denotes the total space of the non-trivial
$S^2$ bundle over $S^2$.
Then $M$ admits infinitely many stable maps into $\bR^2$ with four cusps
of the absolute index one, of the common planar portrait illustrated
in Figure \ref{f:quatro}, such that they are inequivalent 
to each other
up to right-left equivalence.
\end{cor}

\bpf
Let $H_m$ be a Hirzebruch surface with Euler number $m$,
which is a toric surface and 
obtained by pasting four copies of $D^2\times D^2$
along their boundaries by using $\tilde \gamma _i, i=1,\cdots ,4$
with  
$\gamma _1 = m, \gamma _2 = 0, \gamma _3 = -m$,
and $\gamma _4 =0$.
Since $H_m$ is 
the total space of the $S^2$ bundle 
over $S^2$ with Euler number $m$,
and hence 
is diffeomorphic to $S^2\times S^2$ if $m$ is even,
and to $S^2 \tilde\times S^2$ otherwise,
we can obtain a stable map $\tilde f_m$
of $S^2 \times S^2$ or $S^2 \tilde\times S^2$,
by applying the construction in Example \ref{p:toric} for each $m$.
Bellow we show the inequivalence of these maps,
for different positive $m$.

Take a pair of crossing straight lines in the plane
so that it divides the planar portrait into four cusped fans.
Along a line, $M$ is divided 
into two $B_0$ and along the other,
it is divided into $B_m$ and $B_{-m}$,
where $B_k$ denotes the total space of the $D^2$ bundle over $S^2$ with Euler number $k$.
Now take positive integers $m$ and $m'$ 
of the same parity but $m\not=m'$.
If the two maps $\tilde f_m$ and 
$\tilde f_m'$ are right-left equivalent,
then $B_m$ should be mapped by a diffeomorphism of $M$ onto 
one of $B_0,B_{m'} $ or $B_{-m'}$, which is a contradiction. 
\epf 

We can generalise the construction of maps for Hirzebruch surfaces
above as follows.
Let $M$ be the total space of an $S^p$ bundle over $S^q\,
(p,q\geq 1)$ with the structure group $O(p+1)$.
Assume that the bundle has a cross section.
Then by dividing each fibre into a disc neighbourhood of 
the section and the rest, 
and further dividing the base space $S^q$ into two discs,
one can cut $M$ into four pieces of $D^p \times D^q$.
Then in a similar way as before, one can obtain stable maps of $M$
to see that:

\begin{cor}\label{p:ss}
The total space of an $S^p$ bundle over $S^q\,
(p,q\geq 1)$ with the structure group $O(p+1)$
that admits a cross section has the planar portrait
made of four cusped fans as 
drawn in Figure \ref{f:quatro}.
The absolute indices of the cusps are the same and is 
$\max\{p-1,q-1\}$.
\end{cor}

\begin{rem}{\rm
The converse, or a characterisation of such 
sphere bundles by the planar portrait $P$ in Figure \ref{f:quatro}
do not hold, unfortunately.
In actual, one can take a stable map of ${\bf C} P^2 \sharp {\bf C} P^2$
which produces $P$, by taking 
natural connected sum of two projections 
in Example \ref{p:liftp2} and then by performing a cusp ellimination.
On the other hand, the manifold admits no 
such bundle structures, since it is diffeomorphic to neither one
of $S^1 \times S^3, S^1 \tilde\times S^3$ (the non-trivial $S^3$ bundle
 over $S^1$), $S^2\times S^2, S^2\tilde\times S^2$.
}
\erem

\section{Appendix}

{\bf Proof of Lemma \ref{p:Ssmooth}}
By Lemma \ref{p:singset}, $S_f$ is a curve in the $z_0 w_0$ plane.
We may hence reduce the setup to the case of two variables 
$z=z_0,w=w_0$.

Note that $S_f$ passes through the four points $(\pm 1,0)$ and 
$(0,\pm 1)$ on $\bdry \Delta$, at which 
$S_f$ is regular as seen bellow. 
Set $A =\partial _{z_0}S$ and $B=\partial _{w_0}S$,
where 
$\displaystyle\partial _t =\frac{\partial}{\partial t}$.
By calculation,
\begin{eqnarray*}
A 
= & (1-2\ve w)s-(1-2\ve w)z \cdot 2\ve w +2\ve z r. \\
\end{eqnarray*}
Hence 
$A\not= 0$ at $(0,\pm 1)$, and 
by symmetry, $B\not= 0$ at $(\pm 1,0)$.

On the other hand,  $S=0$ is transformed to $ZW=1$
by a coordinate change $(z,w) \mapsto (Z,W)$ of ${\rm Int}\Delta$,
where
$Z=F(z), W=F(w)$, and $\displaystyle F(t)=\frac{(1-2\ve t)t}{\ve (1-t^2)}$.
It is easy to see that one component of $\{S=0\} \cap {\rm Int}\Delta$ 
reaches to $(1,0)$ and $(0,1)$ when $t$ moves toward $\pm 1$,
while the other reaches to $(-1,0)$ and $(0,-1)$.
We hence obtain the lemma.
\epf

{\bf Proof of Lemma \ref{p:1jet}}
One can check the first half easily, 
by calculation on Jacobi matrix $J_f$.
Hence we omit it.

To show that $d^2 f$ is surjective,
note first that $k(z)$ and $k(w)$ are not simultaneously zero on $S_f$.
In actual, suppose $k(z)=k(w)=0$ at a singular point.
Then $q=r=0$, and the Jacobi matrix shows
either $p,p_1, \cdots , p_a =0$ or
$s,s_1,\cdots , s_b =0$. 
They imply
$z=0$ or $w=0$.
But it is a contradiction, since $k(0)=1\not=0$.
Bellow we show the surjectivity in $k(z)\not=0$. 
The proof works in parallel in the other case.

Let $J_f$ be the Jacobi matrix of $f$;
\[
J_f= 2 \cdot\left(
\begin{array}{*8c}
      p & p_1 & \cdots & p_a & q &0 & \cdots & 0 \\
r & 0   & \cdots & 0   & s        & s_1 & \cdots &s_b
\end{array}
\right)\mbox{}_.
\]
Since $q\not=0$, we can change $J_f$
to $J'$ bellow, by a coordinate change.
\[
J'= 2 \cdot\left(
\begin{array}{*8c}
      0 &\multicolumn{2}{c}{ \cdots} & 0         & 1  &0 & \cdots & 0 \\
      S &p_1\,s & \cdots & p_a s   & \ast  & s_1 & \cdots &s_b
\end{array}
\right)\mbox{}_.
\]
Hence at a singular point $P \in S_f$,
$d^2 f$ has the matrix representation
\[
H_P=\left(
\begin{array}{c|*3c|c|*3c}
 A| _P&   &      &   &B| _P&   &      &   \\
\hline
  &C_1&      &   & &   &      &   \\
  &   &\ddots&   & &   &      &   \\
  &   &      &C_a& &   &      &   \\ 
\hline
  &   &      &   & &D_1&      &   \\ 
  &   &      &   & &   &\ddots&   \\ 
  &   &      &   & &   &      &D_b 
\end{array}
\right)\mbox{},
\]
where 
\[
A|_P = ( \partial _{z_0}S)_P, B|_P=(\partial _{w_0}S)_P,
C_i= (\partial _{z_i}p_i s)_P, \mbox{ and }\, D_j = (\partial _{w_j}s_j)_P.
\]

By calculation,
\begin{eqnarray*}
C_1&=&\cdots =C_a=(1-2\ve z_0)(1-2\ve w_0)w_0,  \hspace{1em}\mbox{ and} \\
D_1&=&\cdots =D_b=(1-2\ve z_0). 
\end{eqnarray*}
They do not vanish for small $\ve$
(For $C_i$, note that 
if $w_0=0$ on $S_f$, then $z_0 = \pm1$ by Lemma \ref{p:singset}.
Thus $k(z)=0$. But it is now excluded).
Besides, $A|_P$ and $B|_P$ is not simultaneously zero, 
by the regularity of $S_f$ (Lemma \ref{p:Ssmooth}).
Therefore $d\,{}^2 f$ is a surjection.
\epf

{\bf Proof of Lemma \ref{p:k=0}}
Assume first that $q\not=0$. 
Let $M$ be the matrix of the base change
by which 
$J_f M =J'$.
Strictly,
it is given by;
\[
M=\left(
\begin{array}{c|*3c|c|*3c}
-q        &&&&&&&   \\
\hline
  &-q      &&&&&& \\
  &&\ddots &&&&& \\
  &&& -q     &&&& \\
\hline
p &p_1&\cdots&p_a&\frac{1}{q}&&\\
\hline
&&&&&       1 && \\
&&&&&&  \ddots & \\
&&&&&&& 1
\end{array}
\right).
\]

Now define the basis $\ser{l}{{n-1}}$ of $L=\ker df$ by 
\[ 
<\,l_1,l_2, \cdots, {\partial w}_0,l_{a+2}, \cdots, l_{n-1}>
=<{\partial z}_0,\cdots, {\partial z}_a,
                 {\partial w}_0,\cdots, {\partial w}_b> M,    
\]
where 
$\displaystyle{\partial z}_i=\frac{\partial}{{\partial z}_i}$ and 
$\displaystyle{\partial w}_j=\frac{\partial}{{\partial w}_j}$.
With respect to the basis, 
the map $d^2 f|L$
is represented as 

\[
H_L = H_P \cdot \check M 
    = -
\left(
\begin{array}{*7c}
q A -p B &      -p_1 B & \cdots & - p_a B   &   &      & \\
               &q C_1 &        &           &   &      & \\
               &             & \ddots &           &   &      & \\ 
               &             &        &q C_a&   &      & \\
               &             &        &           &-D_1&      & \\
               &             &        &           &   &\ddots& \\
               &             &        &           &   &      &-D_b
\end{array}
\right),
\]
\noindent where $\check M$ is the matrix obtained 
by removing the $(a+2)$-th column from $M$. 
The $(1,1)$ component is $K$. Hence 
the lemma is obvious if $q\not=0$.
On the other hand, assume that $q=0$, or $k(z) = 0$.
Then $P$ is not in $S_{11}(f)$ nor 
$K(P)=0$, as calculation shows.
\epf

\brem\label{r:absidx} The absolute index of the cusp $P_D$ is ${\rm max}\{a,b\}$, 
since in the representation $H_L$ of $d^2 f |L$ above,
the $a$ diagonals from the $(2,2)$ component and the rest $b$ diagonals
have opposite sign. 
We may assume that $a >b$, by exchanging $z$ and $w$ if necessary.
\erem

{\bf Proof of Lemma \ref{p:Ksingle}}
We reduce the setup to the two variable case 
$z=z_0, w=w_0$ as before, since $S_f$ is in the $z_0w_0$ plane.
Note that $P\in S_{11} (f)$ is not on $\bdry \Delta$,
as seen by checking $K$ at the four points $(\pm 1,0),(0,\pm 1)$ on
$S_f\cap\bdry \Delta $.
Hence by Lemma \ref{p:k=0}, 
we are enough to check the zeros of $K$ in ${\rm Int} \Delta$.

Denote by $S_0$ the connected component of $S_f$ which 
passes through the region $z,w \leq 0$ of $\Delta$,
and $S_1$ the other.
Each component has a unique diagonal point,
since $S_f$ is changed to $ZW=1$ by a 
coordinate change in the proof of Lemma  \ref{p:Ssmooth}
and since the change preserves the diagonal set.
Denote by $P_D$ the diagonal point on $S_1$.
Note that $p^2= q^2$ at $P_D$. But
$p=-q$ is attained 
only at negative $z$, as shown by calculation.
Hence $p=q$ at $P_D$, and thus  
$K$ vanishes at $P_D$.
We show bellow that it is the unique vanishing point of $K$
in ${\rm Int} \Delta$.

We divide $S_1 \cap {\rm Int} \Delta$ 
into the upper and the lower part of $P_D$.

{\sl Claim 1}:\quad 
$0 < p <q$ on the upper part.

{\sl Proof}. 
Note that $p$ is positive, since so is $z$,
and that
$p < p_d=(1-2\ve z)z$, since $z<w$.
It is also easy to see that $p_d < q$ for $0<z<\alpha$,
where we set $P_D=(\alpha, \alpha)$. 
They imply the claim.
\epf

{\sl Claim 2}:\quad 
$0< B<A$ on the upper part, where $A=\bdry_z S$ and $B=\bdry_w S$.

{\sl Proof}.
Note that 
\[ 
(A,B)=(\bdry _Z S, \bdry _W S) 
\left(
\begin{array}{cc}
\bdry _z  F(z) & 0 \\
               0  & \bdry _w F(w)
\end{array}
\right).
\]
The first factor of the right hand side is a gradient vector of $S$, which
is parpendicular
to the curve $ZW=1$, and hence both elements have constant sign,
on this piece.
The diagonal elements in the last factor have also constant sign.
Hence so have both $A$ and $B$.
The sign is positive, as checked
at a regular point $(1,0)$.

On the other hand, we see that
\begin{eqnarray*}
B & = &  (1-2\ve z)(1- 2\ve w)z +2\ve w (\ve z^2-z+\ve)
\end{eqnarray*}
(refer to the equation on $A$ in the proof of Lemma \ref{p:Ssmooth}),
and hence
\begin{eqnarray*}
A-B & = & (w-z)((1-2\ve z)(1-2\ve w)+2\ve^2 zw-2\ve ^2) \\
    & = & (w-z)(6\ve^2 zw-2\ve(z+w)+1-2\ve ^2).
\end{eqnarray*}
The last factor of above is positive on $\hat\Delta$:
it defines a hyperbola 
\[ 
h(z,w):=6(\ve z -\frac{1}{3})(\ve w -\frac{1}{3})+\frac{1}{3}-2\ve ^2 =0
\]
on the
$z$-$w$ plane,
which have no intersection with $\hat\Delta$ for positive small $\ve$,
while $h(0,0)>0$.
Hence we obtain the claim.
\epf

By the claims, $pB<qA$ on the upper part,
and hence $K$ do not vernish.
By symmetry, $K$ has no zeros on $S_1 \cap {\rm Int}\Delta$ 
but for the diagonal point $P_D$.
It is easily checked that $K$ has no zeros on 
${\rm Int}\Delta \cap S_0$,
by using the facts $p<0, q>0$ (since $z,w <0$)
while both $A$ and $B$ are negative.
\epf

{\bf Proof of Lemma \ref{p:KStansvers}}
We can reduce the setup to the two variable case $z=z_0,w=w_0$, as before.
Since 
the tangent direction of
$S_f$ at $P_D$ is $(1,-1)$, 
we are enough to show $\bdry _z K -\bdry _w K \not= 0$ at $P_D$.
Set $J_1=(\bdry _z K) _{P_D}$ and $J_2= (\bdry _w K) _{P_D}$.

Recall that $p=q$ at $P_D$,
$\bdry _z A=\bdry _w B$ at $P_D$
and that $\bdry _w A=\bdry _z B = \bdry _{zw}S$.
By using them,
\begin{eqnarray*}
J_1   
& = & \bdry_z q \cdot A + p \cdot\bdry_z A -\bdry_z p \cdot A -p\cdot\bdry_z B, \quad \mbox{and}\\
J_2 
& = & \bdry_w q \cdot A + p\cdot \bdry_z B -\bdry_w p \cdot A -p\cdot \bdry_z A.
\end{eqnarray*}
Hence 
\begin{eqnarray*} 
J_1 -J_2 & = & (\bdry _z q -\bdry _z p -\bdry _w q + \bdry _w p ) A 
                +2p(\bdry _z A - \bdry_z B) \\
         & = & -(1+2\ve z)A +2p\bdry _z(A - B).
\end{eqnarray*}
On the other hand,
by the equation on $A-B$ in the proof of Lemma \ref{p:Ksingle},
we see that
\begin{eqnarray*} \bdry _z (A-B)_{P_D}=-(6\ve ^2z^2 -4\ve z +1 -2\ve ^2).
\end{eqnarray*}
Substituting this to the equation on $J_1 -J_2$ above, and 
using $A=p(1-2\ve z)$ at $P_D$,
we obtain 
\begin{eqnarray*}
J_1 -J_2 
         =  -z(1-2\ve z)(8\ve ^2 z^2-8\ve z -4\ve ^2 +3).
\end{eqnarray*}

It is now easy to see that  $J_1 -J_2 \not= 0$ at $P_D$,
since the last factor is positive on $|z|\leq 1$ for small $\ve$.
\epf

%\newpage
\vspace{3.0em}

\centerline{ References }
\vspace{1.0em}

\noindent
{[GG]}
       M. Golubitsky and V. Guillemin,
     {\sl Stable mappings and their singularities},
Grad. Texts in Math. 
            {\bf 14},  Springer
            (1973).

\noindent
{[K1]}
       M. Kobayashi,
     {\sl Simplifying certain stable mappings from simply connected
4-manifolds into the plane},
    Tokyo J. Math.
            {\bf 15}
            (1992), 327--349.

\noindent
{[K2]}
      M. Kobayashi,
    {\sl On planar portraits of manifold of depth $2$},
      In preparation.

\noindent
{[K3]}
     M. Kobayashi,
{\sl Twin frames of homotopy spheres and the planar portraits of Morse complexity two},
preprint.

\noindent
{[L1]}
       H. Levine, 
            {\sl   Elimination of cusps},
             Topology 
             {\bf 3}, suppl. 2  
             (1965), 
             263--296.

\noindent
{[L2]}
            H. Levine,
           {\sl   Mappings of manifolds into the plane},
              Amer. J. Math. 
             {\bf 88} 
             (1966), 
             357--365.

\noindent
{[L3]}
            H. Levine,
           {\sl  Classifying immersions into $\bR ^4$ over
stable maps of $3$-manifolds into $\plane$},
              Lect. Notes in Math. 
             {vol. 1157}, Springer,
             1985. 

\noindent
{[Ma]}
            J. Mather,
           {\sl Stability of $C^{\infty}$ mappings: VI.
The nice dimensions}, 
             Proc. Liverpool Singularities - Symposium I, 
             Lect. Notes in Math.
             {vol. 192}, Springer,
            1971.

\noindent
{[Pi]}
          R. Pignoni,
{\sl
On surfaces and their contours
}, Manuscripta Math. 
{\bf 72}(1991), 223--249.

\noindent
{[Th]}
           R. Thom,
{\sl 
Les singularites des applications differentables,
} Annales de l'Institut Fourier,
vol 4 (1955-56), 43-87.

\end{document}